\RequirePackage{fix-cm}
\documentclass[smallextended]{svjour3}       
\smartqed  
\usepackage{graphicx}

\usepackage{amsmath,amscd,amsfonts,amssymb,enumerate, mathrsfs}
\usepackage{color}
\usepackage[colorlinks]{hyperref}

\usepackage{floatrow}
\usepackage{multirow}
\usepackage[table]{xcolor}

\begin{document}
	
	\title{Radius Properties of Harmonic Mappings with Fixed Analytic Part 
		\thanks{The work of Kamajeet Gangania is supported by University Grant Commission, New-Delhi, India  under UGC-Ref. No.:1051/(CSIR-UGC NET JUNE 2017).}
	}
	
	
	\author{      
	         Kamaljeet Gangania 
	}
	
	\authorrunning{Kamaljeet Gangania} 
	
\institute{
		Department of Applied Mathematics, Delhi Technological University,
	Delhi--110042, India \at
	\email{gangania.m1991@gmail.com}           
	\and 
}
	
	\date{Received: date / Accepted: date}

	\maketitle
	
	\begin{abstract}
	We investigate harmonic mappings $f=h+\bar{g}$ defined in the unit disk, where $g$ and $h$ satisfy certain prescribed conditions and the analytic part $h$ belongs to the Ma and Minda class of starlike functions. Certain sharp radius results for univalency, close-to-convexity, fully starlikeness, fully convexity, and radius of strongly starlikeness are established. We also calculate the radius of uniformly starlikeness and convexity for these functions. Several results enhance the well-known radius results. \\
	\end{abstract}
	
	\section{Introduction}
	\label{intro}
	A complex-valued function  $f$ is harmonic if and if $f\in C^2(\Omega)$ (continuous first and second partial derivatives in $\Omega$) with $\Delta{f}=4f_{z\bar{z}}=0$. Suppose that $f$ is a harmonic mapping of a proper simply connected domain $\Omega \subset \mathbb{C}$. Then $f: \Omega \rightarrow \mathbb{C}$ is a harmonic function if and only if the function $f$ has the representation $f=h+\bar{g}$, where $h$ and $g$ are analytic in $\Omega$. Therefore, well-known result of Lewy's~\cite{Lewy-1936} implies that such function $f$ is locally univalent and sense-preserving in $\Omega$ if and only if $|h'(z)|>|g'(z)|$ in $\Omega$. Equivalently, $f$ is sense-preserving local homeomorphism if and only if $J_f(z)>0$, where the jacobian $J_f(z)=|f_z(z)|^2 - |f_{\bar{z}}(z)|^2.$ Further, by Riemann mapping theorem there is a conformal mapping $\varphi$ of $\mathbb{D}$ onto $\Omega$, and it follows that the composition $f\circ \varphi$ maps $\mathbb{D}$ harmonically onto $\Omega$. As a consequence, we can assume that $\Omega$ is the unit disk and that $f$ is sense-preserving in $\mathbb{D}$. Because $J_f>0$ for sense-preserving harmonic mapping $f=h+\bar{g}$, $f_z(0)=h'(0)\neq 0$, and so we arrive at the normalized form $({f(z)-f(0)})/{f_z(0)}$ of $f$.
	Thus, the analytic and co-analytic parts of the harmonic mapping $f=h+\bar{g}$ in $\mathbb{D}$ can be written, respectively,
	\begin{equation}\label{seriesexpansion}
	h(z)=z+\sum_{n=2}^{\infty}a_nz^n \quad \text{and} \quad g(z)=\sum_{n=1}^{\infty}b_nz^n.
	\end{equation}
	
	Denote by $\mathcal{S_{H}}$ (see~\cite{Clunie-Sheil-Small-1984}) the class of all complex-valued harmonic univalent and sense-preserving mappings $f$ in the unit disk $\mathbb{D}$ with $f(0)=0=f_z(0)=1$. Let $\mathcal{A}$ denotes the class of analytic functions in $\mathbb{D}$ and normalized by the condition $h(0)=0=h'(0)-1$. We note that $\mathcal{S_{H}}$ reduces to $\mathcal{S}$, the class of normalized univalent analytic functions in $\mathbb{D}$ whenever the co-analytic part of $f$ is zero. Furthermore, for $f=h+\bar{g} \in \mathcal{S_H}$ with $g'(0)=b$ and $|b|<1$ (because $J_f(0)=1-|b|^2>0$), the function 
	$F=({f-\overline{bf}})({1-|b|^2}) \in \mathcal{S_H}.$
	The univalency of $F$ is obtained by applying an affine mapping to $f$. Observe that $F_{\bar{z}}(0)=0$.  Therefore, we may restrict our attention to the subclass
	\begin{equation*}
	\mathcal{S}_{\mathcal{H}}^{0} = \{f \in \mathcal{S_H}: f_{\bar{z}}(0)=0 \}.
	\end{equation*}
	The condition $f_{\bar{z}}(0)=0$ is equivalent to the assertion that the second complex dilatation $\omega_f(z)= \overline{f_{z}(z)}/ f_{z}(z)$ of $f$ is zero. Clearly,
	$\mathcal{S}\subset \mathcal{S}_{\mathcal{H}}^{0} \subset \mathcal{S_H}$. Although both the classes $\mathcal{S_H}$ and $\mathcal{S}_{\mathcal{H}}^{0}$ are known to be normal families, only $\mathcal{S}_{H}^{0}$ is compact with respect to the topology of locally uniform convergence. The geometric subclasses of $\mathcal{S_H}$ consisting of convex, starlike and close-to-convex functions in $\mathbb{D}$ are denoted $\mathcal{C_H}$, $\mathcal{S}_{\mathcal{H}}^{*}$ and $\mathcal{K_H}$ respectively. Let $\mathcal{C}_{\mathcal{H}}^{0}$, $\mathcal{S}_{\mathcal{H}}^{*0}$ and $\mathcal{K}_{\mathcal{H}}^{0}$ denote the subclasses of $\mathcal{C_H}$, $\mathcal{S}_{\mathcal{H}}^{*}$ and $\mathcal{K_H}$ with the condition $f_{\bar{z}}(0)=0$ respectively. 
		
	In 1992, Ma and Minda~\cite{minda94} introduced and studied the subclasses of $\mathcal{S}$ of starlike and convex univalent functions :
	\begin{equation*}\label{mindasclass}
	\mathcal{S}^{*}(\psi)=\left\{h\in\mathcal{A}: \frac{zh'(z)}{h(z)}\prec\psi(z)\right\}\; \text{and}\; 	\mathcal{C}(\psi)=\left\{h\in\mathcal{A}: 1+\frac{zh''(z)}{h'(z)}\prec\psi(z)\right\},
	\end{equation*}
	where the univalent function $\psi$ defined as $\psi(z)=1+p_1z+p_2z^2+\cdots$ has positive real in $\mathbb{D}$. The symbol $ \prec $ denotes the usual subordination. Note that $\mathcal{S}^{*}(({1+z})({1-z}))$ and $\mathcal{C}(({1+z})({1-z}))$ respectively reduce to the well-known classes $\mathcal{S}^{*}$ and $\mathcal{C}$ of starlike and convex functions. For some current radius results, see~\cite{ganga_cmft2021,harmonic-IJST2022,kamal-mediter2022,ganga-KQsi-2021,janow,Kumar-cardioid,sinefun,sokol09} and their references.
	
	For $f\in \mathcal{S}^*(\psi)$ or $\mathcal{C}(\psi)$, the starlikeness or convexity of the range set $f(\mathbb{D})$ is inherited by the level sets $f(\mathbb{D}_r)$ for all $0<r<1$. However, the hereditary property of starlikeness and convexity for harmonic function does not hold. For this, Chuaqui et al.~\cite{Chuaqui} defined the notion of fully starlike and fully convex harmonic mappings.
	\begin{definition}\cite{Chuaqui}
	 A harmonic function $f=h+\bar{g}$ with $h$ and $g$ of the form \eqref{seriesexpansion} is said to be fully starlike of order $0\leq \beta<1$ if each circle $C_r=\{z: |z|=r \}$ is one-to-one mapped onto a curve that bounds a starlike domain with respect to the zero and satisfy
	\begin{equation*}
	\frac{\partial}{\partial \theta} (\arg{f(re^{i\theta})}) >\beta; \quad 0\leq \theta<2\pi, \;0<r<1.
	\end{equation*}
	Similary, $f$ is said to be fully convex of order $\beta$, if $f$ maps one-to-one each circle $C_r$ onto a curve that bounds a convex domain and satisfies
	\begin{equation*}
	\frac{\partial}{\partial \theta}\{\arg( \frac{\partial}{\partial \theta} f(re^{i\theta})) \} >\beta; \quad 0\leq \theta<2\pi, \;0<r<1.
	\end{equation*}
	\end{definition}
In $1969$, Bernardi generalized Alexander's theorem by introducing the function $\phi_{\nu}: \mathbb{D} \rightarrow \mathbb{C}$ defined by
\begin{equation*}
\phi_{\nu}(z)= \frac{\nu+1}{z^{\nu}} \int_{0}^{z} \xi^{\nu-1}\phi(\xi)d\xi,
\end{equation*}
where $\phi$ is analytic in $\mathbb{D}$, with $\phi(0)=0=\phi'(0)-1$. According to the Rad\'{o}-Kaneser-Choquet theorem (see~\cite{Clunie-Sheil-Small-1984}), fully convex functions of order $\beta$ are necessarily univalent in $\mathbb{D}$. However, fully harmonic starlike functions need not be univalent in $\mathbb{D}$. The well-known Alexander's relation $\phi \in \mathcal{S}^*$ if and only if  $\int_{0}^{z}\phi(\xi)/\xi d\xi \in \mathcal{C}$ between the classes of starlike and convex function. Sheil-Small~\cite{Sheil-Small-1990} extended this result to the harmonic mappings and proved that
\begin{theorem}\cite{Sheil-Small-1990}(Sheil-Small)\label{Sheil-Small-Alexander}
	If $f=h+\bar{g}: \mathbb{D}\rightarrow \mathbb{C}$ fixes zero, is univalent, and has a starlike range, and $H$ and $G$ are the analytic functions in $\mathbb{D}$ defined by $zH'(z)=h(z),$ $zG'(z)=-g(z)$, $H(0)=G(0)=0$, then $F=H+\overline{G}$ is univalent, and has a convex range. 
\end{theorem}
For $\nu\geq 0$, $\mu \in\mathbb{D}$ and $f\in \mathcal{S_H}$, Muir~\cite{Muir-2017} generalized the Bernardi integral operator to the harmonic functions, by defining $\bigwedge_{\nu,\mu}[f]: \mathbb{D} \rightarrow \mathcal{C}$ as follows:
\begin{equation*}
\bigwedge_{\nu,\mu}[f](z)= \frac{\nu+1}{z^{\nu}} \int_{0}^{z} \xi^{\nu-1}h(\xi)d\xi  + \mu \overline{  \frac{\nu+1}{z^{\nu}} \int_{0}^{z} \xi^{\nu-1}g(\xi)d\xi }.
\end{equation*}
Clearly, $\bigwedge_{0,-1}$ is the harmonic analogue of Alexander's operator, given in Theorem~\ref{Sheil-Small-Alexander}. In fact, the classes $\mathcal{C_H}$ and $\mathcal{S}_{\mathcal{H}}^{*}$ are preserved under the operator $\bigwedge_{0,-1}$. However, $\mathcal{C_H}$ and $\mathcal{S}_{\mathcal{H}}^{*}$ are not necessarily preserved under $\bigwedge_{0,1}$. Indeed, $\bigwedge_{0,1}[K] \not\in \mathcal{S}_{\mathcal{H}}^{*0}$, where the harmonic Koebe function $K(z)$ is given by
\begin{equation*}
K(z)=\frac{z-\frac{1}{2}z^2+\frac{1}{6}z^3}{(1-z)^3} + \overline{ \frac{\frac{1}{2}z^2+\frac{1}{6}z^3}{(1-z)^3} }.
\end{equation*}

There has been an interplay between the sense-preserving harmonic mappings and their analytic part. Mocanu~\cite{Mocanu-2011 } conjectured that if $h$ and $g$ are analytic functions in $\mathbb{D}$ that satisfy $g'(z)=zh'(z)$ and $\Re(1+zh''(z)/h'(z))>-1/2$, then the harmonic function $f=h+\bar{g}$ is univalent in $\mathbb{D}$.  Bshouty and Lyssaik~\cite{Bshouty-Lyzzaik-2011} proved that a sense-preserving harmonic function $f=h+\bar{g}$ where $h$ and $g$ of the form \eqref{seriesexpansion} is necessarily univalent in $\mathbb{D}$ if the dilatation of $f$ is $w_f(z)=z$ and $h$ satisfies $\Re(1+zh''(z)/h'(z))>-1/2$ for all $z\in \mathbb{D}$.

Recently, harmonic functions were studied with different dilatation and analytic part $h$ belonging to different geometric subclasses, see~\cite{Bharendhar-Ponnu-2014,Bshouty-Lyzzaik-2011,Ponnu-Kaliraj2015,prajapat-Mathi2022,Ankur-Sumit-2022}. Recall that a domain $\Omega$ is called close-to-convex if the complement of $\Omega$ can be written as a union of half-lines such that corresponding open half-lines are disjoint. A function $f\in \mathcal{S_{H}}$ which maps $\mathbb{D}$ onto a close-to-convex domain is called a {\it close-to-convex function}. Clunie and Sheil-Small~\cite{Clunie-Sheil-Small-1984} proved that
\begin{theorem}\cite{Clunie-Sheil-Small-1984}\label{Sheil}
	Let $f=h+\bar{g}$ be a locally univalent harmonic function in $\mathbb{D}$ and $h+\epsilon g$ $(|\epsilon|\leq 1)$ is convex in $\mathbb{D}$. Then $f$ is univalent and close-to-convex in $\mathbb{D}$.
\end{theorem}  

Note that, in particular, if $h$ is starlike in Theorem~\ref{Sheil}, then $f$ need not be univalent in $\mathbb{D}$. The Koebe function $k(z)=z/(1-z)^2$ is starlike. However, the harmonic function given by
\begin{equation*}
 F(z)=\frac{z}{(1-z)^2} +\overline{\frac{\ln(1-z)}{2} +\frac{z(1-2z)}{2(1-z)^2}}
\end{equation*}
is sense-preserving in $\mathbb{D}$ as $|g'(z)/h'(z)|=|z(1+z)/2|<1$ but not univalent in $\mathbb{D}$. Further, the harmonic function
\begin{equation*}
G(z)= \frac{z}{(1-z)^2} +\overline{\frac{z}{2(1-z)} }
\end{equation*}
is sense-reversing in $\mathbb{D}$ as $|g'(-3/4)/h'(-3/4)|=7/2>1$, but $G$ is not univalent in $\mathbb{D}$ and also $g(z)=\frac{1}{2}(1-z)h(z)$. Hence, it is natural to determine the radius of univalency of harmonic functions which are not necessarily sense-preserving when $h\in \mathcal{S}^*(\psi)$ and $g(z)=\varphi(z) h(z)$, where $\varphi\in\mathcal{B}_s$, the class of analytic functions $\varphi: \mathbb{D}\rightarrow \mathbb{D}$ such that $|\phi(z)|\leq 1$ for all $z\in \mathbb{D}$.

We investigate the harmonic functions and the operator $\bigwedge_{0,1}$ for the radius of close-to-convexity and radius of univalency when the analytic part is a Ma-Minda starlike function (that is, $h\in \mathcal{S}^*(\psi)$) when the analytic and co-analytic part satisfy certain conditions. Further, we derive the hereditary radius of fully starlikeness, fully convexity, strongly starlikeness and strongly convexity of such harmonic functions. Unless otherwise stated, we assume that the function $\psi$ has positive coefficients in its power series expansion throughout the discussion.

\section{ Harmonic mappings with Analytic part from the class of Ma-Minda starlike functions}

The radius of univalency and the radius of close-to-convexity of complex-valued harmonic functions $f=h+\bar{g}$, where $h$ and $g$ are given by \eqref{seriesexpansion} satisfying the coefficient conditions
\begin{equation*}
 |a_n|\leq \frac{(2n+1)(n+1)}{6} \quad \text{and} \quad |b_n|\leq \frac{(2n-1)(n-1)}{6}, \quad n\geq 2
\end{equation*}
is found to be $0.112903$, see~\cite{Kalaj-Poonu-Vuorinen-2014}. Recently, Ponnusamy and Kaliraj~\cite[Theorem~3.8]{Ponnu-Kaliraj2015} proved that for sense-preserving harmonic functions $f=h+\bar{g}$ such that $h$ is univalent in $\mathbb{D}$, where $h$ and $g$ are given by \eqref{seriesexpansion}, the radius of univalency and close-to-convexity equals to $2-\sqrt{3} \approx 0.267949$.

Now, we obtain the radius of sense-preserving and univalency for the harmonic functions $f=h+\bar{g}$ such that $g(z)=\varphi(z)h(z)$ and $h \in \mathcal{S}^*(\psi)$, where $\varphi \in \mathcal{B}_s$. For this, we need to recall
\begin{lemma}\cite{prajapat-Mathi2022}\label{sense-preservingLemma}
	 Let $r\in(0,1)$ and  $f=h+\bar{g}$ be a sense-preserving harmonic mapping in $\mathbb{D}_r=\{z: |z|<r  \}$. Further let $h$ be univalent convex in $\mathbb{D}$. Then $f$ is univalent in $\mathbb{D}_r$.
\end{lemma}	

\begin{theorem}\label{majorization-use}
	Let $f=h+\bar{g}$ be the harmonic mapping in $\mathbb{D}$ such that $h \in \mathcal{S}^*(\psi)$ and $g(z)=\varphi(z)h(z)$, where $\varphi \in \mathcal{B}_s$. Further, assume that
	$$\min_{|z|=r}|\psi(z)|=m(r).$$  
	Then $f$ is sense-preserving and univalent in $\mathbb{D}_R$, where $R=\min\{r_{\psi}, r_{c} \}$ and $r_{\psi}$ is the unique postive root of the equation
	\begin{equation*}
	(1-r^2)m(r) -2r=0
	\end{equation*}
	and $r_{c}$ is the radius of convexity of functions $h$ in $\mathcal{S}^*(\psi)$. If $r_c\geq r_{\psi}$, then the result is sharp for the case $m(r)=\psi(-r)$.
\end{theorem}
\begin{proof}
	Since $h$ is convex (univalent) in $|z|<r_c$. From Lemma~\ref{sense-preservingLemma}, it suffices to prove that $f$ is sense-preserving in $|z|<R$. As $h\in \mathcal{S}^*(\psi)$, by definition we have
	\begin{equation*}
	\frac{zh'(z)}{h(z)}=\psi(\omega(z)), \quad \omega\in \mathcal{B}_s,
	\end{equation*}
	Since $\psi$ is univalent, it follows by Maximum-minimum principle of modulus with the hypothesis $\min_{|z|=r}|\psi(z)|=m(r)$ that
	\begin{equation}\label{starlike-bound}
	\left|\frac{h(z)}{h'(z)} \right| \leq \frac{|z|}{\min_{|z|=r}|\psi(z)|}= \frac{r}{m(r)}.
	\end{equation}  
	Since $g(z)=\varphi(z)h(z)$ implies that
	\begin{equation}\label{maj-inequality}
	 |g'(z)|= |h'(z)| \left\{|\varphi(z)| + \left|\frac{h(z)}{h'(z)}\right||\varphi'(z)| \right\}.
	\end{equation}
	Now by employing the Schwarz-Pick inequality
	\begin{equation*}
	 |\varphi'(z)| \leq \frac{1-|\varphi(z)|^2}{1-|z|^2}
	\end{equation*}
	in the Eq.~\eqref{maj-inequality}, and using \eqref{starlike-bound} gives that
 $    |g'(z)| \leq T(r) |h'(z)|,$
	where 
	\begin{equation*}
	T(r):=\left\{|\varphi(z)| + \frac{r}{m(r)} \left(\frac{1-|\varphi(z)|^2}{1-|z|^2} \right) \right\}.
	\end{equation*}
	Now using \cite[Theorem~5]{ganga_cmft2021} we have $|T(r)|<1$ such that
	\begin{equation*}
	 |g'(z)| < |h'(z)| \quad \text{in}\; |z|< r_\psi,
	\end{equation*}
	where the radius $r_{\psi}$ is sharp and is the smallest root of the equation 
	$$(1-r^2)m(r) -2r=0$$ in $(0,1)$ This in view of Lemma~\ref{sense-preservingLemma}, shows that $f$ is sense-preserving and univalent in $|z|<R=\min\{r_c, r_{\psi}\}$. Further, if $r_c> r_{\psi}$ and $m(r)=\psi(-r)$, then \cite[Proof of Theorem~2]{ganga_cmft2021} implies that
	\begin{equation*}
	|g'(z)| > |h'(z)| \quad \text{for all}\; z=r> r_\psi,
	\end{equation*}
	which completes the proof of sharpness. \qed
\end{proof}	

Now, for different choices of the function $\psi$, we expolre Theorem~\ref{majorization-use}. Note that for harmonic functions $f=h+\bar{g}$ such that $g$ is majorized by $h$ and $h \in \mathcal{S}^*(\psi)$, we obtain the improved radius of univalency in Corollary~\ref{all-result}. For sense-preserving harmonic functions with univalent analytic part, radius of univalency is $2-\sqrt{3}> 3-\sqrt{8}$ the conjectured radius of convexity of $f\in \mathcal{S}_{H}$, see \cite{Sheil-Small-1990}. Let us recall an known class of starlike functions defined by $$\mathcal{S}^*[M]:=\{ h \in \mathcal{A}: \left|\frac{zf'(z)}{f(z)}-M \right|<M \}.$$ Such functions implies that $1+\Re\frac{zh''(z)}{h'(z)} \geq \frac{(1+m)r}{(1-r)(1+mr)} $, where $m=1-\frac{1}{M}.$
\begin{corollary}\label{all-result}
Let $f=h+\bar{g}$ be the harmonic mapping in $\mathbb{D}$ such that $h \in \mathcal{S}^*(\psi)$ and $g(z)=\varphi(z)h(z)$, where $\varphi \in \mathcal{B}_s$. 
Then $f$ is sense-preserving and univalent in $\mathbb{D}_R$, where $R=\min\{r_{\psi}, r_{c} \}$ and $r_{\psi}$ is the unique postive root of the equation
\begin{equation*}
(1-r^2)m(r) -2r=0
\end{equation*}
and $r_{c}$ is the radius of convexity of functions $h$ in $\mathcal{S}^*(\psi)$. The result follows for each one of the following cases:
	\begin{itemize}
		\item [$(i)$]  	$R=2-\sqrt{3}=r_{\psi}=r_c\approx 0.267949$ when $\psi(z)= (1+z)/(1-z)$.
	
		\item [$(ii)$] 	$R=r_{\psi}\approx 0.3524$ is the root of $(1-r)^3(1+r)^2-4r^2=0$ for $\psi(z)=\sqrt{1+z}$.

		\item [$(iii)$] $R= r_{\psi}\approx 0.3237$ is the root of $(1-r^2)-2re^{r}=0$ for $\psi(z)=e^z$.

		\item [$(iv)$] $R=r_{\psi}\approx 0.358473$ is the root of $(1-r^2)(\sqrt{1+r^2}-r)-2r$ when $\psi(z)=z+\sqrt{1+z^2}$.

		\item [$(v)$]  $R=r_{\psi}\approx 0.358473$ is the root of $(1-r^2)-r(1+e^r)=0$ for $\psi(z)=2/(1+e^{-z})$.
		
		\item [$(vi)$] \label{sin} $R=r_c\approx 0.345$ for $\psi(z)=1+\sin{z}$. (see, \cite[Corollary~2.4]{sinefun})
		
		\item [$(vii)$] $R=r_{\psi}$ is the root of $(1+r)(1-r)^2-2r(1+mr)=0$ when $h\in \mathcal{S}^*[M]$.
	\end{itemize}

\end{corollary}

Next to study the radius of fully starlikeness and convexity for functions in $\mathcal{S}_{\mathcal{H}}^{0}$ with the condition $g(z)=\varphi(z)h(z)$, we invoke a result of Jahangiri.
\begin{lemma}\cite{Jahangiri-suffConvxHar1998,Jahangiri-suffStarHar1999} \label{suff-FulStrCnxHarmonic}
	Let the harmonic function $f$ be given by $f(z)=h(z)+\overline{g(z)}= z+\sum_{n=2}^{\infty}a_n z^n+ \overline{\sum_{n=2}^{\infty}b_nz^n}$ in $\mathbb{D}$ and $0\leq \beta<1$.
	\begin{enumerate}
		\item [$(i)$] The function $f$ is fully starlike of order $\beta$, if 
		$$\sum_{n=2}^{\infty}(n-\beta)|a_n| +\sum_{n=2}^{\infty}(n+\beta)|b_n|\leq 1-\beta.$$
		\item [$(ii)$] The function $f$ is fully convex of order $\beta$, if
		$$\sum_{n=2}^{\infty}n(n-\beta)|a_n| +\sum_{n=2}^{\infty}n(n+\beta)|b_n|\leq 1-\beta.$$
	\end{enumerate}
\end{lemma}	

\begin{remark}
	In the case when the coefficient of power series of the function $\psi$ in the class $\mathcal{S}^*(\psi)$ are negative, then without loss of generality, we can replace the function $h_\psi(z)=z \exp\int_{0}^{z}(\psi(t)-1)/tdt:=z+\sum_{n=2}^{\infty}a_nz^n$ by $\hat{h}_\psi(z):=z+\sum_{n=2}^{\infty}|a_n|z^n$ and same for its derivative, in what follows.
\end{remark}	
\begin{theorem}\label{fullystarlike}
	Let $f=h+\bar{g}$ be a  harmonic function such that $h\in \mathcal{S}^*(\psi)$ in $\mathbb{D}$ and $g(z)=\varphi(z)h(z)$, where $\varphi\in \mathcal{B}_s$, $h$ and $g$ are given by \eqref{seriesexpansion}. Then for $\beta\in[0,1)$, $f$ is fully starlike of order $\beta$ in $|z|<R_\beta =\min\{1/3, r_{\psi} \}$, where $r_{\psi}$ is the smallest postive root of the equation
	\begin{equation*}
	4h'_{\psi}(r)+\beta=5
	\end{equation*}
	in $(0,1),$ where $h_{\psi}(z)= z\exp \left(\int_{0}^{z}(\psi(t)-1)/{t}dt \right)$.
\end{theorem}
\begin{proof}
	Let $0<r<1$. Consider the function
	\begin{equation*}
	 f_r(z)=\frac{1}{r} f(rz)= z+ \sum_{n=2}^{\infty}a_n r^{n-1}z^n +\overline{\sum_{n=2}^{\infty}b_nr^{n-1}z^n}.
	\end{equation*}
	Now, it suffices to prove that the function $f_r$ is fully starlike of order $\beta$ in $|z| \leq R_{\beta}$. Since $g(z)=\varphi(z) h(z)$, where $\varphi\in \mathcal{B}_s$, it follows from \cite[Corollary~1]{kamal-mediter2022} that in $|z|\leq 1/3$
	\begin{equation}\label{bound1/3}
	|b_n| \leq |a_n|.
	\end{equation}
	Invoking the result \cite[Corollary~1]{kamal-mediter2022} for the subordination $h(z)/z \prec h_{\psi}(z)/z$ (see~\cite{minda94}) gives together with inequality~\eqref{bound1/3} that for the disk $|z|\leq 1/3$
	 \begin{equation}\label{final-bound1/3}
	 |a_n| \leq |\hat{a_n}|,
	 \end{equation}
	 where $h_{\psi}(z)= z\exp \left(\int_{0}^{z}(\psi(t)-1)/{t}dt \right)= z+\sum_{n=2}^{\infty}\hat{a_n}z^n$. Therefore, from \eqref{bound1/3} and \eqref{final-bound1/3}
	 \begin{align*}
	 S_1 &= \sum_{n=2}^{\infty}\left( \frac{n-\beta}{1-\beta} \right) |a_n|  r^{n-1} + \sum_{n=2}^{\infty}\left( \frac{n-\beta}{1-\beta} \right) |b_n|  r^{n-1}\\
	  & \leq \sum_{n=2}^{\infty}\left( \frac{n-\beta}{1-\beta} \right) |\hat{a_n}|  r^{n-1} + \sum_{n=2}^{\infty}\left( \frac{n-\beta}{1-\beta} \right) |\hat{a_n}|  r^{n-1}\\
	  &= \frac{4}{1-\beta} \sum_{n=2}^{\infty}n|\hat{a_n}|r^{n-1}
	  = \frac{4}{1-\beta} (h'_{\psi}(r)-1)
	 \end{align*}
	 holds in $|z|\leq 1/3$.  Hence, by Lemma~\ref{suff-FulStrCnxHarmonic}, we see that $S_1\leq 1$ holds in $|z|\leq R_{\beta}=\min\{ 1/3, r_{\psi}\}$, where $r_{\psi}$ is the smallest positive root of the equation
	 \begin{equation*}
	 \frac{4}{1-\beta} (h'_{\psi}(r)-1)=1,
	 \end{equation*}
	which implies that $f$ is fully starlike of order $\beta$ in $|z|\leq R_{\beta}$. \qed	 
\end{proof}	
\begin{corollary}
	Let $f$ satisfies the hypothesis of Theorem~\ref{fullystarlike}.  Then for $\beta\in[0,1)$, $f$ is fully starlike of order $\beta$ in $|z|<R_\beta =\min\{1/3, r_{\psi} \}$, where $r_{\psi}$ is the smallest postive root of the equation:
	\begin{itemize}
		\item [$(i)$]  	$4h'_{\psi}(r)+\beta=5$ for $\psi(z)= \frac{1+Dz}{1+Ez}$, where $-1\leq E<D\leq1$ and
		\begin{equation*}
		h'_{\psi}(r)=
		\left\{
		\begin{array}
		{lr}
		\psi(r)(1+E r)^{\frac{D-E}{E}}, & E\neq 0; \\
		
		\psi(r)e^{Dr}, & E=0.
		\end{array}
		\right.
		\end{equation*}

		\item [$(ii)$] $4\psi(r) e^{(\int_{0}^{r}\frac{\psi(t)-1}{t}dt)} +(\beta-5)=0$ for $\psi(z)=(\frac{1+z}{1-z})^{\eta}$, where $0<\eta \leq1$.
		
		\item [$(iii)$] 	$(16\sqrt{1+r}) e^{2(\sqrt{1+r}-1)} +(\beta-5){(\sqrt{1+r}+1)^2}=0$ for $\psi(z)=\sqrt{1+z}$ and $r_{\psi}<1/3$ for all $\beta$.

		\item [$(iv)$] $4(1+re^r) e^{e^r -1}+ (\beta-5)=0$ for $\psi(z)=1+ze^z$ and $r_{\psi}<1/3$ for all $\beta$.

		\item [$(v)$] $8 \psi(r)e^{(\psi(r) -1)} +(\beta-5)(\sqrt{1+r^2}+1)=0$ for $\psi(z)=z+\sqrt{1+z^2}$ and $r_{\psi}<1/3$ for all $\beta$.
		
	\end{itemize}
	
\end{corollary}

\begin{figure}[h]
	\begin{tabular}{c}
		\includegraphics[scale=0.4]{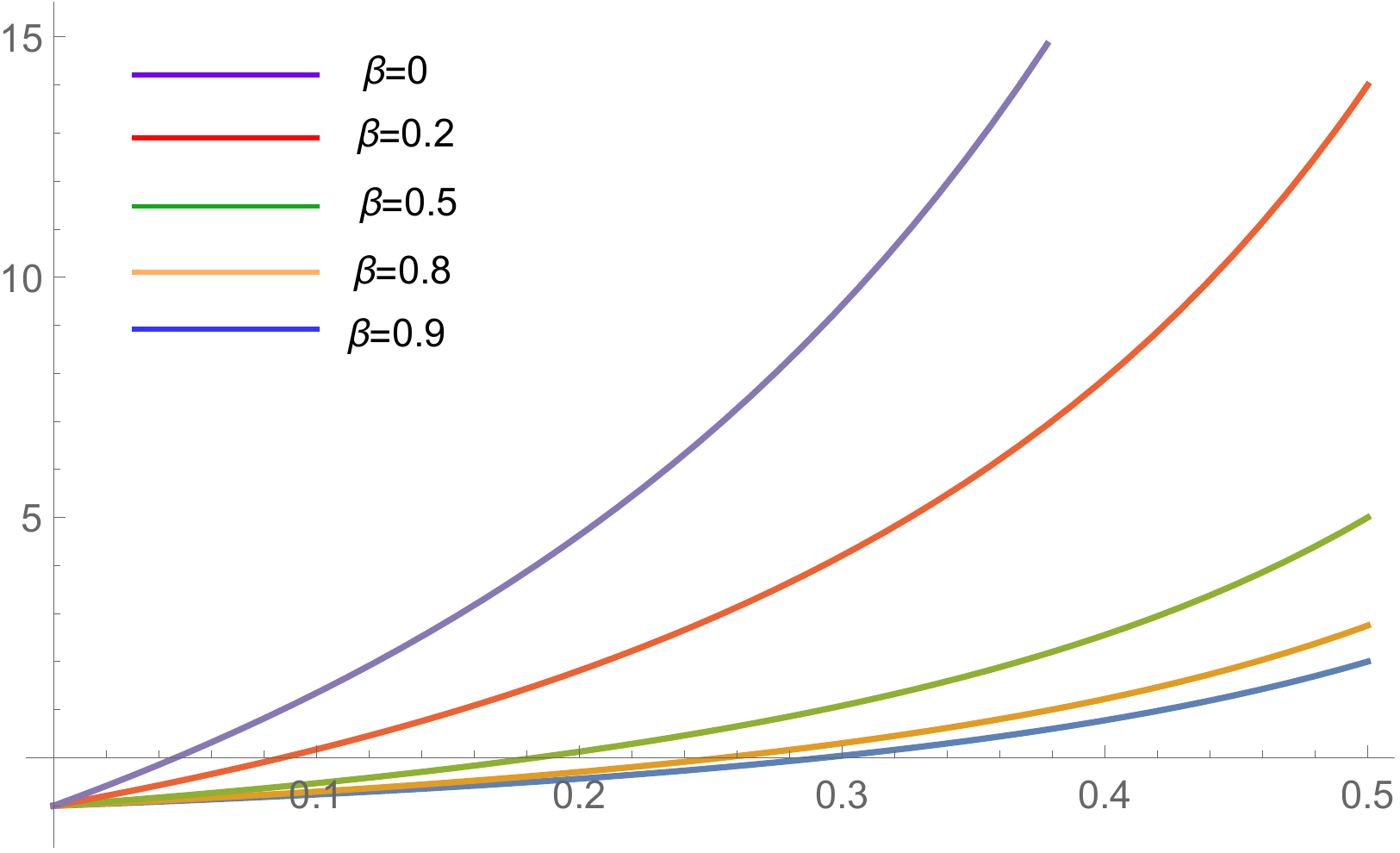}
	\end{tabular}
	\caption{Graph for the roots of equation for various choices of $\beta$ with $\psi(z)=(1+z)/(1-z)$.}\label{CorFFullyStar}
\end{figure}

\begin{corollary}\label{Ffullystarlike}
	Let $f$ satisfies the hypothesis of Theorem~\ref{fullystarlike}. Then for $\beta\in[0,1)$, $F$ given by
	\begin{equation*}
	F(z)= \bigwedge_{0,1}[f](z) =\int_{0}^{z}\frac{h(\xi)}{\xi}d\xi +\overline{\int_{0}^{z}\frac{g(\xi)}{\xi}d\xi  }
	\end{equation*}
	is fully starlike of order $\beta$ in $|z|<R_\beta =\min\{1/3, r_{\psi} \}$, where $r_{\psi}$ is the smallest postive root of the equation
	\begin{equation*}
	\frac{2}{1-\beta}\left (\frac{h_{\psi}(r)-r}{r} \right) =1
	\end{equation*}
	in $(0,1),$ where $h_{\psi}(z)= z\exp \left(\int_{0}^{z}(\psi(t)-1)/{t}dt \right)$. See, Figure~\ref{CorFFullyStar}.
\end{corollary}

 Hernandez and Martin~\cite{Harnandez} (also see, \cite[Theorem~2.7]{Nagpal-ravichandran-2012}) proved that a sense-preserving harmonic function $f=h+\bar{g}$ is fully starlike in $\mathbb{D}$ if the analytic functions $h+\epsilon g$ are starlike in $\mathbb{D}$ for each $|\epsilon|=1$. Now combining this result with a majorization result~\cite[Theorem~2]{ganga_cmft2021} of Gangania and Kumar, we get that

\begin{theorem}\label{improvedfullystarlike}
	Let $h$ and $g$ are given by \eqref{seriesexpansion} and $f=h+\bar{g}$ be a  harmonic function such that $h\in \mathcal{S}^*(\psi)$ in $\mathbb{D}$ and $g(z)=\varphi(z)h(z)$, where $\varphi\in \mathcal{B}_s$. Further let $m(r)=\min_{|z|=r}\Re\psi(z)$. Then for $\beta\in[0,1)$, $f$ is fully starlike of order $\beta$ in $|z|<R_\beta =\min\{ r_1, r_2 \}$, where $r_1$ and $r_2$, respectively are the smallest positive roots of the following equations
\begin{equation*}
m(r)-\frac{r}{1-r}=\beta
\end{equation*}
and
\begin{equation*}
(1-r^2)m(r) -2r=0
\end{equation*}
in $(0,1),$ where $h_{\psi}(z)= z\exp \left(\int_{0}^{z}(\psi(t)-1)/{t}dt \right)$.		
\end{theorem} 

\begin{figure}[h]
	\begin{tabular}{c}
		\includegraphics[scale=0.22]{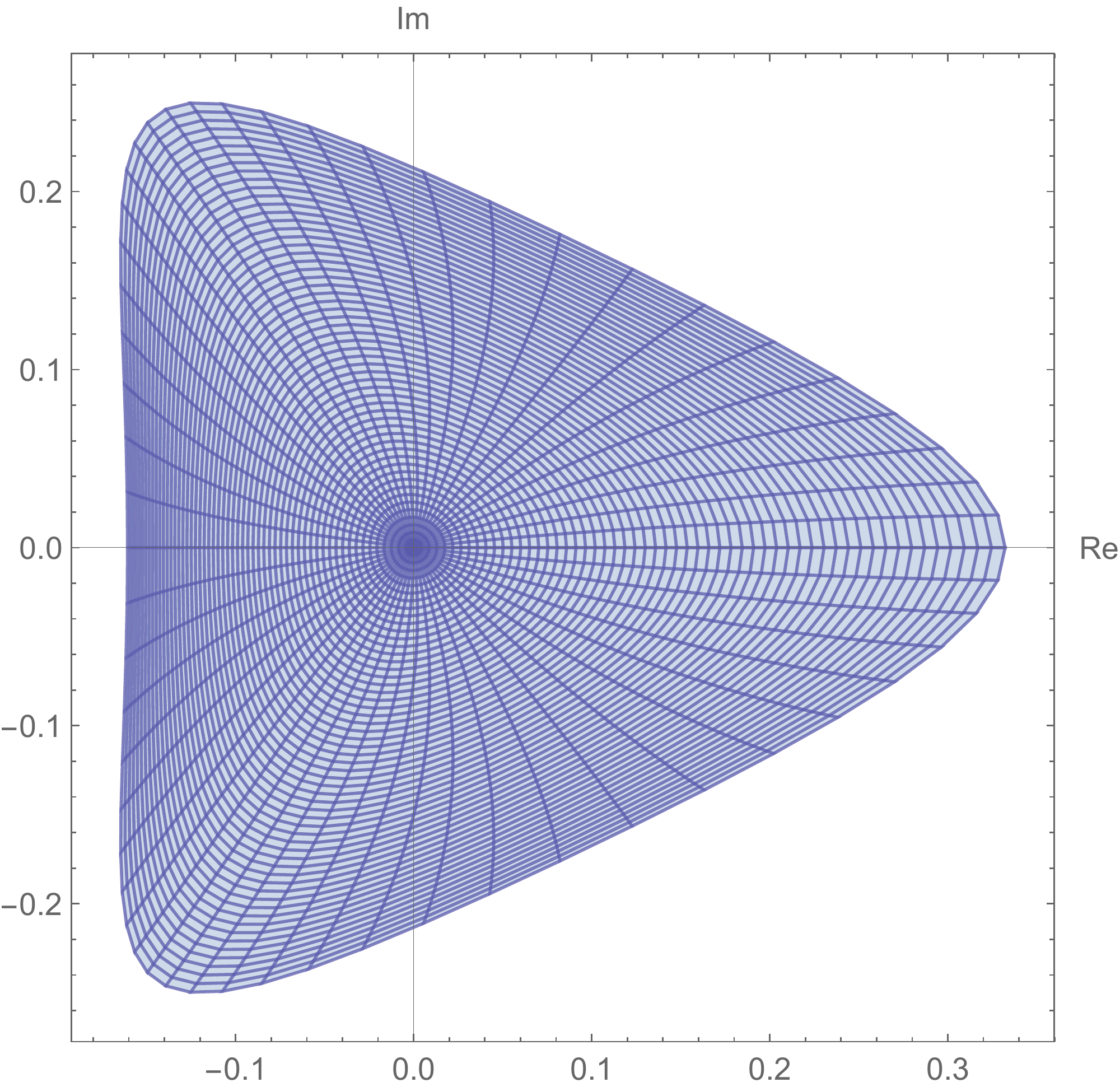}
	\end{tabular}
	\begin{tabular}{l}
		\includegraphics[scale=0.22]{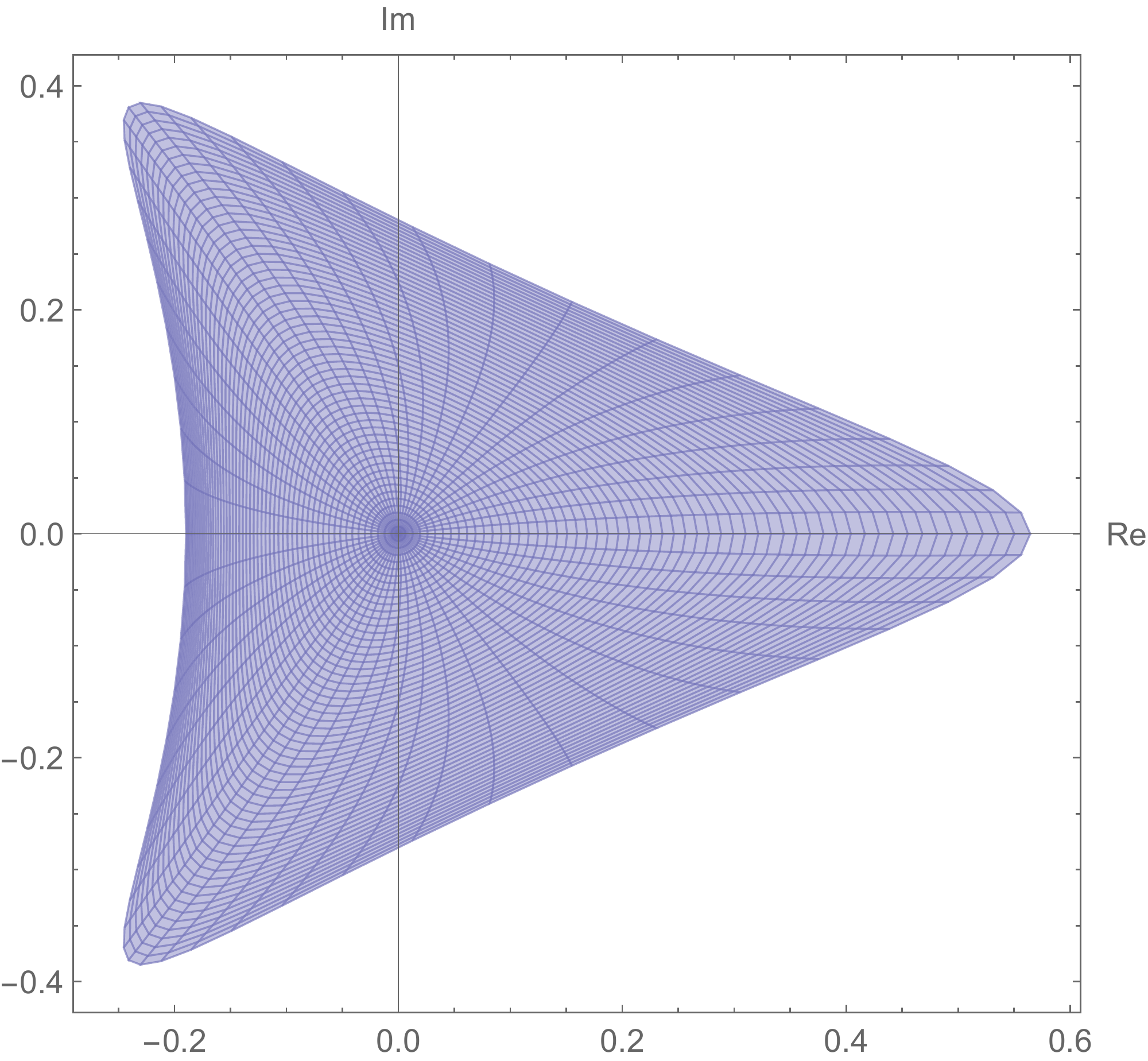}
	\end{tabular}
	\caption{Images of Harmonic function $f_0(z)+\overline{zf_0(z)}$ with $f_0(z)=\frac{4z\exp(2\sqrt{1+z}-2)}{(1+\sqrt{1+z})^2}$.}\label{lemiscateharmonic}
\end{figure}

In Figure~\ref{lemiscateharmonic}, the left side graph is for $|z|< 0.238778$ in View of Theorem~\ref{fullystarlike} and right side is for $|z|<0.3524$ in view of Theorem~\ref{improvedfullystarlike} for the case when $\psi(z)=\sqrt{1+z}$.

In Theorem~\ref{improvedfullystarlike} if we include the condition that $f$ is sense-preserving harmonic function then the radius can be further improved, which is stated in the following result.
\begin{corollary}
	Let $f=h+\bar{g}$ be a sense-preserving harmonic function such that $h\in \mathcal{S}^*(\psi)$ in $\mathbb{D}$. Further let $m(r)=\min_{|z|=r}\Re\psi(z)$. Then for $\beta\in[0,1)$, $f$ is fully starlike of order $\beta$ in $|z|<R_\beta$, which is the smallest positive root of the following equation
	\begin{equation*}
	m(r)-\frac{r}{1-r}=\beta.
	\end{equation*}	
\end{corollary}	
For the choice of $\psi(z)=\sqrt{1+z}$ in above corollary, gives $R_0=0.43016> r_2=0.3524$, where $r_2$ is defined in Theorem~\ref{improvedfullystarlike}.

\begin{theorem}\label{fullyconvex}
	Let $f=h+\bar{g}$ be a harmonic function such that $h\in \mathcal{S}^*(\psi)$ in $\mathbb{D}$ and $g(z)=\varphi(z)h(z)$, where $\varphi\in \mathcal{B}_s$, $h$ and $g$ are given by \eqref{seriesexpansion}. Then for $\beta\in[0,1)$, $f$ is fully convex of order $\beta$ in $|z|<C_\beta =\{1/3, r_{\psi} \}$, where $r_{\psi}$ is the smallest postive root of the equation
	\begin{equation*}
	G'(r)* h'_{\psi}(r)= \frac{5-\beta}{4}
	\end{equation*}
	in $(0,1)$, where $G(z)=z/(1-z)$ and the symbol $``*"$ represents the convolution.
\end{theorem}
\begin{proof}
 Similar to the proof of Theorem~\ref{fullystarlike}, it is enough to show that $f_r$ is fully convex in $|z|\leq C_{\beta}$. Now using the result \cite[Corollary~1]{kamal-mediter2022} on $g(z)=\varphi(z) h(z)$ and \eqref{final-bound1/3}, we see that
 	\begin{align*}
 	S_2 &= \sum_{n=2}^{\infty}\left( \frac{n(n-\beta)}{1-\beta} \right) |a_n|  r^{n-1} + \sum_{n=2}^{\infty}\left( \frac{n(n+\beta)}{1-\beta} \right) |b_n|  r^{n-1}\\
 	& \leq \sum_{n=2}^{\infty}\left( \frac{n(n-\beta)}{1-\beta} \right) |\hat{a_n}|  r^{n-1} + \sum_{n=2}^{\infty}\left( \frac{n(n+\beta)}{1-\beta} \right) |\hat{a_n}|  r^{n-1}\\
 	&= \frac{4}{1-\beta} \sum_{n=2}^{\infty}n^2|\hat{a_n}|r^{n-1}
  	\end{align*}
  	holds in $|z|\leq 1/3.$ Note that using convolution, we can re-write
  	\begin{equation*}
  	\sum_{n=2}^{\infty}n^2|\hat{a_n}|r^{n-1}= 	(G'(r)-1)* (h'_{\psi}(r)-1),
  	\end{equation*}
  	which implies that $S_2\leq 1$ holds in $|z|\leq C_{\beta}=\min\{1/3, r_{\psi}\}$, where $r_{\psi}$ is the least positive root in $(0,1)$ of the equation
  	\begin{equation*}
  	(G'(r)-1)* (h'_{\psi}(r)-1)= \frac{1-\beta}{4},
  	\end{equation*}
  	 where $G(z)=z/(1-z)$. Hence, Lemma~\ref{suff-FulStrCnxHarmonic} says that $f$ is fully convex of order $\beta$ in $|z|\leq C_{\beta}$. \qed
\end{proof}	

\begin{corollary}\label{Ffullyconvex}
	Let $f$ satisfies the hypothesis of Theorem~\ref{fullyconvex}. Then for $\beta\in[0,1)$, $F$ given by
	\begin{equation*}
	F(z)= \bigwedge_{0,1}[f](z) =\int_{0}^{z}\frac{h(\xi)}{\xi}d\xi +\overline{\int_{0}^{z}\frac{g(\xi)}{\xi}d\xi  }
	\end{equation*}
	is fully convex of order $\beta$ in $|z|<C_\beta =\min\{1/3, r_{\psi} \}$, where $r_{\psi}$ is the smallest postive root of the equation
	\begin{equation*}
	2h'_{\psi}(r)+\beta=3
	\end{equation*}
	in $(0,1),$ where $h_{\psi}(z)= z\exp \left(\int_{0}^{z}(\psi(t)-1)/{t}dt \right)$.
\end{corollary}

\begin{corollary}\label{Cases-fullyconvex}
	Let $f$ satisfies the conditions of Theorem~\ref{fullyconvex}. Then for $\beta\in[0,1)$, $f$ is fully convex of order $\beta$ in $|z|<C_\beta =\{1/3, r_{\psi} \}$, where $r_{\psi}$ is the smallest positive root of the equation
	\begin{equation*}
	G'(r)* h'_{\psi}(r)= \frac{5-\beta}{4}
	\end{equation*}
	in $(0,1)$, where $G(z)=z/(1-z)$ and
	\begin{itemize}
		\item [$(i)$]  for $\psi(z)= \frac{1+Dz}{1+Ez}$, where $-1\leq E<D\leq1$:
		\begin{equation*}
		h'_{\psi}(r)=
		\left\{
		\begin{array}
		{lr}
		\psi(r)(1+E r)^{\frac{D-E}{E}}, & E\neq 0; \\
		
		\psi(r)e^{Dr}, & E=0.
		\end{array}
		\right.
		\end{equation*}
		
		\item [$(ii)$]  for $\psi(z)=\sqrt{1+z}$: \quad 	$h'_{\psi}(r)=\frac{4\sqrt{1+r}}{(\sqrt{1+r}+1)^2} e^{2(\sqrt{1+r}-1)}$.

		\item [$(iii)$] for $\psi(z)=1+ze^z$: \quad $h'_{\psi}(r)=(1+re^r)e^{e^r -1}$ .

		\item [$(iv)$] for $\psi(z)=z+\sqrt{1+z^2}$: \quad $h'_{\psi}(r)=\frac{2(r+\sqrt{1+r^2})}{\sqrt{1+r^2}+1} e^{(r+\sqrt{1+r^2}-1)}$.
		
	\end{itemize}
	
\end{corollary}

\begin{figure}[h]
	\begin{tabular}{c}
		\includegraphics[scale=0.25]{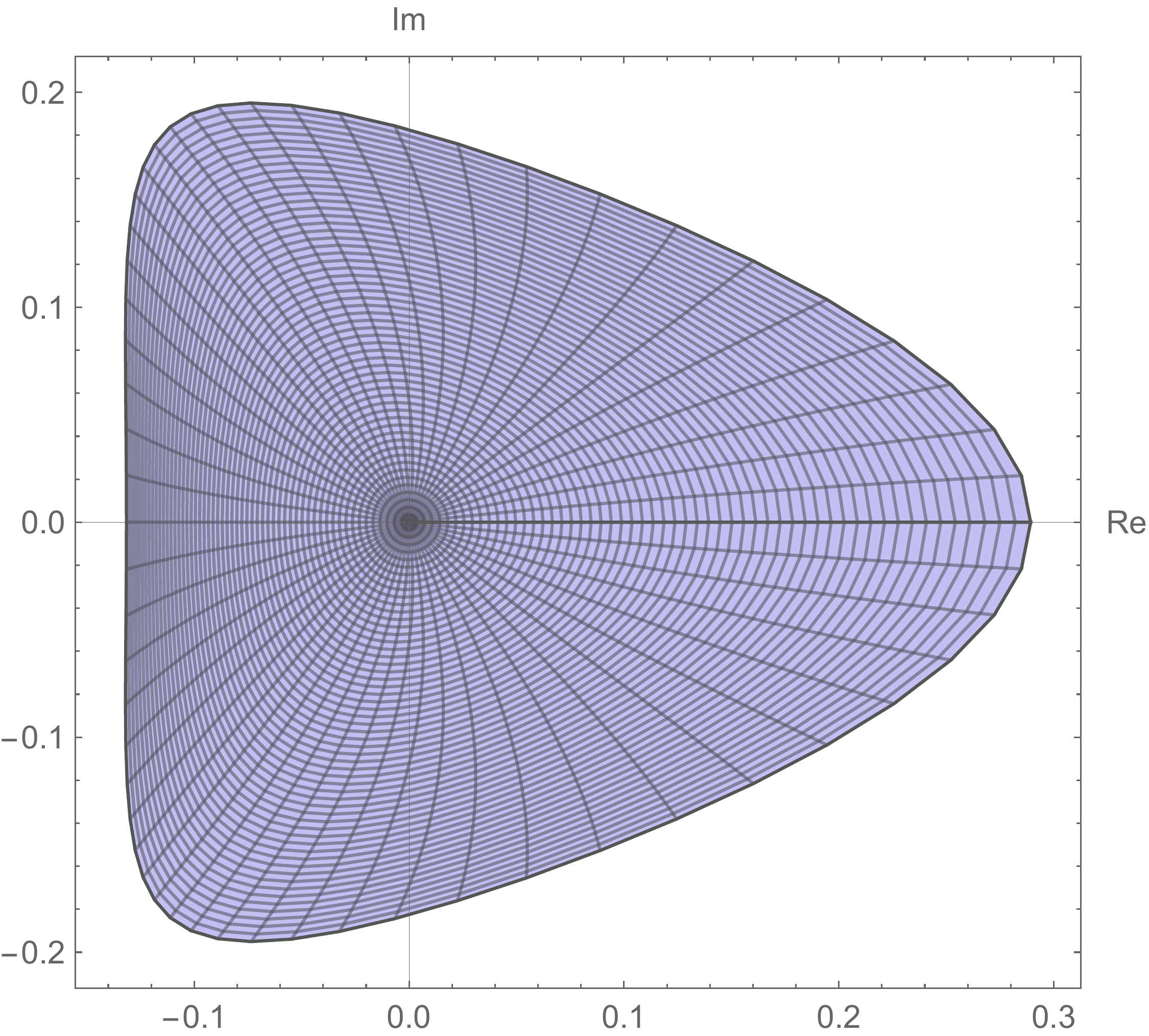}
	\end{tabular}
	\caption{Image of the harmonic function $z\exp(e^z-1)+ \overline{z^2\exp(e^z-1)}$ in $|z|\leq 0.195106$.}\label{sharpfullyconvexity}
\end{figure}

The following table show that the radius $C_\beta =r_{\psi}<1/3$ in part $(iii)$ of Corollary~\ref{Cases-fullyconvex} is sharp, see Figure~\ref{sharpfullyconvexity}. 
\begin{table}[h]
	\begin{tabular}{|p{1.5cm}|p{2cm}|p{2cm}|p{2cm}|p{2cm}|}
		\hline
		\textbf{}     & $S_5$ & $S_{10}$ & $S_{20}$ & $S_{n\rightarrow \infty}$\\ \hline
		$\beta=0$    & $0.1952$ & $0.195106$ & $0.195106$ & $0.195106$ \\ \hline
		$\beta=1/2$   & $0.181806$ & $0.181742$ & $0.181742$ & $0.181742$ \\ \hline
		$\beta=0.9$   & $0.17048$ & $0.170435$ & $0.170435$ & $0.170435$ \\ \hline
		             &  &   &  &$r_{\psi}<1/3$ \\
		\hline	
	\end{tabular}
	\caption{Behaviour of the roots of the equation $G'(r)* h'_{\psi}(r)= \frac{5-\beta}{4}$ when $\psi(z)=1+ze^z$ in Corollary~\ref{Cases-fullyconvex}.}
\end{table}

For the sense-preserving harmonic functions $f=h+\bar{g}$ such that $f_{\bar{z}}(0)=0$ and $h\in \mathcal{K}(c)$, where $$\mathcal{K}(c):= \left\{1+\frac{zh''(z)}{h'(z)}>c,\; -\frac{1}{2}\leq c<1 \right\}$$
is the class of convex functions of order $c$ (such functions are univalent and convex in one direction in $\mathbb{D}$), Ponnusamy and Kaliraj~\cite[Theorem~4.6]{Ponnu-Kaliraj2015} proved that $f$ is univalent and close-to-convex in $|z|<r_c$, where $r_c=1/3$ for $-1/2\leq c\leq c_0\approx-0.391827$, and $r_c=\cos{\pi c}$ for $c_0\leq c \leq 0$, where $c_0$ is the root of $\cos{\pi c}=1/3$.

Let us now consider the function $\psi(z)=\sqrt{1+z}$. Then $\mathcal{S}^*(\psi)$ reduces to the class of lemniscate starlike functions, denoted by $\mathcal{SL}^*$.
\begin{example}\label{leminsacate-harmonic}
	Let $h \in \mathcal{SL}^*$ and $\varphi \in \mathcal{B}_s$. If $g'(z)=\varphi(z)h'(z)$ for all $z\in \mathbb{D}$, then the harmonic mappings $f=h+\bar{g}$ are univalent and close-to-convex in $|z|< r_0\approx0.3119$ is the unique positive roof of $1-6r+12r^2-11r^3+4r^4=0$
	in $(0,1)$.
\end{example}
\begin{proof}
	Let $Q_{\delta}(z)=h(z)-\delta g(z)$, where $|\delta|\leq 1$. Then we have
	\begin{equation*}
	 Q'_{\delta}(z)=(1-\delta \varphi(z))h'(z),
	\end{equation*}
	which further gives
	\begin{equation}\label{convex-equation}
	1+\frac{zQ''_{\delta}(z)}{Q'_{\delta}(z)}= 1+\frac{zh''(z)}{h'(z)}-\frac{\delta z\varphi'(z)}{1-\delta \varphi(z)}.
	\end{equation}
	Let 
	$ p(z)={zh'(z)}/{h(z)},$
	then using Lemma~\cite{sokol09}, we get the sharp inequality
	\begin{equation}\label{lemsicate}
	\Re\left( 1+\frac{zh''(z)}{h'(z)} \right) \geq \sqrt{1-r} -\frac{r}{\sqrt{1-r}}. 
	\end{equation}
	Using inequality~\eqref{lemsicate} in \eqref{convex-equation}, we get 
	\begin{equation}\label{final-convex}
	\Re\left( 1+\frac{zQ''_{\delta}(z)}{Q'_{\delta}(z)} \right) \geq \sqrt{1-r} -\frac{r}{\sqrt{1-r}} -\Re \left( \frac{\delta z\varphi'(z)}{1-\delta \varphi(z)} \right). 
	\end{equation}
	The Schwarz-Pick inequality says that $\varphi\in \mathcal{B}_s$ satisfies
	\begin{equation*}
	 |\varphi'(z)| \leq \frac{1-|\varphi(z)|^2}{1-|z|^2}, \quad\text{which gives}\quad
	 \left|\frac{\delta \varphi'(z)}{1-\delta \varphi(z)} \right| \leq \frac{|z|}{1-|z|}.
	\end{equation*}
 Hence, 
	\begin{equation*}
	\Re\left(\frac{\delta \varphi'(z)}{1-\delta \varphi(z)} \right) \leq \sqrt{1-r} -\frac{r}{\sqrt{1-r}}
	\end{equation*}
	for $|z|=r< r_0$, where $r_0\approx0.3119$ is the smallest positive roof of $1-6r+12r^2-11r^3+4r^4=0$. Hence, in view of \eqref{final-convex}, $f$ is univalent and close-to-convex in $|z|<r_0$. \qed
\end{proof}	

If the dilatation $\varphi=g'/h'$ of harmonic mapping $f=h+\bar{g}$ has the form $\varphi(z)=e^{i\theta} z^n$ $(\theta \in \mathbb{R}), n\in \mathbb{N}$, then we get
\begin{example}
	Let $f=h+\bar{g}$ be the harmonic mapping in $\mathbb{D}$ such that $h \in \mathcal{SL}^*$ and $g'(z)=e^{i\theta} z^n h'(z)$. Then $f$ is univalent and close-to-convex in $|z|<r_n$, where $r_n$ is the smallest positive root of the equations
	\begin{equation*}
	(1-4r+4r^2) -r^n(2-8r+8r^2)+r^{2n}(1-n^2+(n^2-4)r+4r^2)=0
	\end{equation*} 
	in $(0,1)$.
\end{example}

We now generalize Example~\ref{leminsacate-harmonic}. In fact, the following result can also be compared with \cite[Theorem~3.8]{Ponnu-Kaliraj2015}. 
For certain results for the harmonic function $f=h+\bar{g}$ to be univalent and close-to-convex, we refer to see \cite{Bharendhar-Ponnu-2014,Bshouty-Joshi-2013,Bshouty-Lyzzaik-2011}.
\begin{theorem}\label{Gen-leminsacate-harmonic}
Let $h \in \mathcal{S}^*(\psi)$ and $\varphi \in \mathcal{B}_s$. If $g'(z)=\varphi(z)h'(z)$ for all $z\in \mathbb{D}$, then the harmonic functions $f=h+\bar{g}$ are univalent and close-to-convex in $|z|< r=\min\{r_0, r_c\}$, where $r_c$ is radius of convexity for the class $\mathcal{S}^*(\psi)$ and $r_0$ is the minimal positive root of the equation
$$(1-r)H_{\psi}(r)-r=0$$
with
$$  H_{\psi}(r):=	\min_{0<|z|=r}\Re\left( 1+\frac{zh''(z)}{h'(z)} \right).$$

\end{theorem}
\begin{proof}
	In Theorem~\ref{leminsacate-harmonic}, we notice that one requires the sharp inequality~\eqref{lemsicate} in case when $h\in \mathcal{S}^*(\psi)$ which deals with the finding of radius of convexity for functions $h$. So without loss of generality, we may assume a continuous function $H_{\psi}(r)$ such that there exists some $r_c\in (0,1)$ in such a way that $H_{\psi}(r)\geq0$ for $r\in[0,r_c]$ and $H_{\psi}(r)<0$ for $r\in (r_c, r_c+\epsilon )$ for some $\epsilon>0$, and the following sharp inequality holds
	\begin{equation*}
	\Re\left( 1+\frac{zh''(z)}{h'(z)} \right) \geq H_{\psi}(r)> 0, \quad\text{for all}\quad |z|=r\in [0,r_c),
	\end{equation*} 
    where $H_{\psi}(0)=1$. Hence, further following the footsteps of the Proof of Example~\ref{leminsacate-harmonic}, we get the desired result. \qed
\end{proof}	

\begin{corollary}
	Let $h \in \mathcal{S}^*_{SG}$ in Theorem~\ref{Gen-leminsacate-harmonic}. If $g'(z)=\varphi(z)h'(z)$ for all $z\in \mathbb{D}$, where $\varphi \in \mathcal{B}_s$ then the harmonic mappings $f=h+\bar{g}$ are univalent and close-to-convex in $|z|< r_0 \approx 0.3729$, where $r_0$ is the unique positive root of the equation
	$$\frac{2-re^r}{1+e^r}-\frac{r}{1-r}=0.$$
	
\end{corollary}	

Let us now recall a basic definition to proceed further.
\begin{definition}
A locally univalent  function $f=h+\bar{g}$ is said to be uniformly starlike in $\mathbb{D}$, if $f$ is fully starlike in $\mathbb{D}$, and maps every circular arc $\Gamma_\zeta$ contained in $\mathbb{D}$ with center $\zeta\in \mathbb{D}$ onto the arc $f(\Gamma_\zeta)$ which is starlike with respect to $\zeta$. Similarly, we have definition of uniformly convex function $f=h+\bar{g}$ in $\mathbb{D}$. 
\end{definition}

Recall the following important sufficient conditions for functions to be uniformly starlike and convex:
\begin{lemma}\cite{Ponnu-Prajapat-Sairam-2015} \label{suff-uniStrCnxHarmonic}
	Let the harmonic function $f=h+\bar{g} \in \mathcal{S}_{\mathcal{H}}^{0}$. Then the following holds:
	\begin{enumerate}
		\item [$(i)$] The function $f$ is uniformly starlike, if 
		$\sum_{n=2}^{\infty} n(|a_n|+|b_n|) \leq 1/2.$
		\item [$(ii)$] The function $f$ is uniformly convex, if
		$\sum_{n=2}^{\infty} n(2n-1)(|a_n|+|b_n|) \leq 1.$
	\end{enumerate}
\end{lemma}	

Using the Lemma~\ref{suff-uniStrCnxHarmonic}, for harmonic functions $f \in \mathcal{S}_{\mathcal{H}}^{0}$, radius problems for the operator $\bigwedge_{0,1}[f]$ were studied in \cite{Gosh-Vasudeva-2019} with the prescribed conditions on the coefficient bounds $|a_n|$ and $|b_n|$ for the analytic and co-analytic part. 
We now prove the result when $g$ is majorized by some function $h$ in $\mathcal{S}^*(\psi)$.
\begin{theorem}\label{uniformlyconvex}
	Let $f=h+\bar{g} \in \mathcal{S}_{\mathcal{H}}^{0}$ such that $g(z)=\varphi(z)h(z)$ and $h\in \mathcal{S}^*(\psi)$ in $\mathbb{D}$, where $\varphi\in \mathcal{B}_s$. Then $f$ is uniformly convex on the disk  $|z|<r_{uc} =\min\{1/3, r_{\psi} \}$, where $r_{\psi}$ is the smallest positive root of the equation
	\begin{equation*}
	2(G'(r)-1)*(h'_{\psi}(r)-1)-(h'_\psi(r)-1)=\frac{1}{2}
	\end{equation*}
	in $(0,1),$ where $G(z)=z/(1-z)$ and $h_{\psi}(z)= z\exp \left(\int_{0}^{z}(\psi(t)-1)/{t}dt \right)$.
\end{theorem}
\begin{proof}
	For $0<r<1$, let 
	\begin{equation*}
	f_r(z)=r^{-1}(rz)= z+\sum_{n=2}^{\infty}a_nr^{n-1}z^n +\overline{ \sum_{n=2}^{\infty} b_n r^{n-1}z^n}, \quad (z\in \mathbb{D}).
	\end{equation*}
	In view of Lemma~\ref{suff-uniStrCnxHarmonic}, let us consider the summation
	\begin{equation*}
	S= \sum_{n=2}^{\infty}|a_n|r^{n-1} +\sum_{n=2}^{\infty} |b_n| r^{n-1}.
	\end{equation*}
	Thus, from \eqref{bound1/3} and \eqref{final-bound1/3}, it follows that
	\begin{align*}
	S \leq \sum_{n=2}^{\infty} 2n(2n-1)|\hat{a_n}| r^{n-1}
	&=: 2[2(G'(r)-1)*(h'_{\psi}(r)-1) -(h'_{\psi}(r)-1)]\\
	&\leq 1
	\end{align*}
	holds in $|z|=r\leq \min\{ 1/3, r_{\psi}\}$, where $r_{\psi}$ is the least positive root of the equation
	\begin{equation*}
	2(G'(r)-1)*(h'_{\psi}(r)-1)-(h'_\psi(r)-1)=\frac{1}{2}.
	\end{equation*}
	This proves that $f$ is uniformly convex on the disk $|z|< r_{uc}$. \qed
\end{proof}	

\begin{theorem}\label{uniformlystarlike}
	 Let $f$ satisfies the hypothesis of Theorem~\ref{uniformlyconvex}. Then $f$ is uniformly starlike on the disk  $|z|<r_{us} =\min\{1/3, r_{\psi} \}$, where $r_{\psi}$ is the smallest positive root of the equation
	\begin{equation*}
	h'_{\psi}(r)-\frac{5}{4}=0
	\end{equation*}
	in $(0,1),$ where $h_{\psi}(z)= z\exp \left(\int_{0}^{z}(\psi(t)-1)/{t}dt \right)$.
\end{theorem}
\begin{proof}
	It follows using the technique of Theorem~\ref{uniformlyconvex} with inequalities~\eqref{bound1/3}, \eqref{final-bound1/3} and the application of  Lemma~\ref{suff-uniStrCnxHarmonic}. \qed
\end{proof}	

\begin{corollary}\label{F-unistarlike}
	Let $f$ satisfies the hypothesis of Theorem~\ref{uniformlyconvex}. Then $F$ given by
	\begin{equation*}
	F(z)= \bigwedge_{0,1}[f](z) =\int_{0}^{z}\frac{h(\xi)}{\xi}d\xi +\overline{\int_{0}^{z}\frac{g(\xi)}{\xi}d\xi  }
	\end{equation*}
	is uniformly starlike on the disk  $|z|<r_{us} =\min\{1/3, r_{\psi} \}$, where $r_{\psi}$ is the smallest positive root of the equation
	\begin{equation*}
	4h_{\psi}(r)-5r=0
	\end{equation*}
	in $(0,1),$ where $h_{\psi}(z)= z\exp \left(\int_{0}^{z}(\psi(t)-1)/{t}dt \right)$.
\end{corollary}

The above corollary yields that
\begin{corollary}
	Let $f$ satisfies the hypothesis of Theorem~\ref{uniformlyconvex} with $\psi(z)=1+ze^z$. Then $F$ given by
	\begin{equation*}
	F(z)= \bigwedge_{0,1}[f](z) =\int_{0}^{z}\frac{h(\xi)}{\xi}d\xi +\overline{\int_{0}^{z}\frac{g(\xi)}{\xi}d\xi  }
	\end{equation*}
	is uniformly starlike on the disk  $|z|<r_{us} \approx 0.201424$ is the unique positive root of the equation
	\begin{equation*}
	4 r\exp(e^r -1) -5r=0.
	\end{equation*}
\end{corollary}

In $2022$, the notion of hereditarily strongly starlike harmonic functions appeared.
\begin{definition}\cite{Ma-Ponnu-Sugawa}
	Let $f$ be the complex-valued harmonic function $f=h+\bar{g}$ with $h$ and $g$ of the form \eqref{seriesexpansion}. The function $f$ is called {it hereditarly strongly starlike of order $\alpha \in (0,1)$ } if it is sense-preserving and univalent on $\mathbb{D}$ and if $f(\mathbb{D}_r)$ is a strongly starlike domain of order $\alpha$ for each $0<r<1$.   
\end{definition}	

For such functions sharp bounds for $|a_n|$ and $|b_n|$ are yet to be known. We determine the radius of hereditarliy strongly starlikeness of order $\alpha$. 
\begin{theorem}\label{stronglyharmonic}
	Let the harmonic functions $f=h+\bar{g} \in \mathcal{S}_{\mathcal{H}}^{0}$ such that $g(z)=\varphi(z)h(z)$ and $h\in \mathcal{S}^*(\psi)$ in $\mathbb{D}$, where $\varphi\in \mathcal{B}_s$. Then $f\in \mathcal{SS}_{\mathcal{H}}^{0}(\alpha) $ on the disk  $|z|<r_{ss} =\min\{1/3, r_{\psi}(\alpha) \}$, where $r_{\psi}(\alpha)$ is the smallest positive root of the equation
	\begin{equation*}
    \left( \frac{h_{\psi}(r)-r}{r} \right)* M_{\alpha}(r)-2\sin{\frac{\pi \alpha}{2}}=0,
	\end{equation*}
	where $0<\alpha <1$, 
	and $M_{\alpha}(r)= \sum_{n=2}^{\infty}(A_n(\alpha)+B_n(\alpha))r^{n-1}$ with
	$A_n(\alpha)=n-1+|n-e^{-i \pi \alpha}| $ and $B_n(\alpha)=n+1+|n+e^{i\pi \alpha}|.$ 
\end{theorem}
\begin{proof}
	 The result follows by invoking \cite[Theorem~4.2]{Ma-Ponnu-Sugawa}, and proceeding similar to the proof of Theorem~\ref{uniformlyconvex}. \qed
\end{proof}	

\section*{Conflict of interest}
	The authors declare that they have no conflict of interest.

	

\end{document}